\documentclass{amsart}

\usepackage{amsthm}
\newtheorem{thm}{Theorem}[section]

\theoremstyle{definition}

\usepackage{amssymb}
\usepackage{amsmath}
\usepackage{amsfonts}

\usepackage{verbatim}

\newcommand\inprod[2]{\left(#1,#2\right)}
\renewcommand\Im{\mathrm{Im}}

\title[Rank-$40r$ Extremal Even Unimodular Lattices]{Configurations of Rank-$40r$ Extremal Even Unimodular Lattices ($r=1,2,3$)}
\author{Scott Duke Kominers and Zachary Abel}

\address{Department of Mathematics, Harvard University}
\email{SDK: \texttt{kominers@fas.harvard.edu}\phantom{e: }ZA: \texttt{zabel@fas.harvard.edu}}


\date{December 5, 2007}

\begin{document}
\begin{abstract}
  We show that if $L$ is an extremal even unimodular lattice of rank
  $40r$ with $r=1,2,3$ then $L$ is generated by its vectors of norms $4r$ and $\mbox{4r+2}$.
  Our result is an extension of Ozeki's analogous result for the case $r=1$.
\end{abstract}

\maketitle

\section{Introduction}
A \textit{lattice} of \textit{rank} $n$ is a free $\mathbb{Z}$-module
of rank $n$ equipped with a positive-definite inner product
$\inprod{\cdot}{\cdot}: L\times L\to\mathbb{R}$.  The \textit{dual} of
$L$, denoted $L^*$, is the set
\begin{equation*}L^* = \left\{y\in L\otimes\mathbb{R}: \forall x\in L, \inprod{x}{y}\in\mathbb{Z} \right\} ,\end{equation*}
which itself forms a lattice of the same rank as $L$. For a lattice
vector $x\in L$, we call $\inprod{x}{x}$ the \textit{norm} of $x$. A
lattice $L$ is \textit{integral} if $\inprod{x}{x'}\in\mathbb{Z}$ for
all $x,x'\in L$, i.e. if and only if $L\subseteq L^*$. An integral
lattice is said to be \textit{unimodular} if it is self-dual
($L=L^*$).

A lattice $L$ is called \textit{even} if and only if every lattice
vector has an even integer norm, i.e. $\inprod{x}{x}\in 2\mathbb{Z}$
for $x\in L$.  An even lattice is automatically integral by the
familiar parallelogram identity,
$2\inprod{x}{x'}=\inprod{x+x'}{x+x'}-\inprod{x}{x}-\inprod{x'}{x'}$.

Lattices that are simultaneously even and unimodular are especially
rare.  Indeed, such a lattice's rank must be divisible by $8$. Sloane
proved that if $L$ is an even unimodular lattice of rank $n$ then
the minimal (nonzero) norm in $L$ is bounded by
\begin{equation}
  \min_{\stackrel{x\in L}{x\neq 0}}\inprod{x}{x}
  \leq  2\lfloor  n/24\rfloor+2
  \label{minnorm}
\end{equation}
(see \cite[p. 194, Cor. 21]{Conway:SPLAG}). An even unimodular lattice
of rank $n$ is called \textit{extremal} if it attains the bound
\eqref{minnorm}.

Ozeki \cite{Ozeki1, Ozeki3} showed that if $L$ is an extremal even
unimodular lattice of rank $32$ or $48$ then $L$ is generated by its
vectors of minimal norm.  The first author \cite{SK} showed analogous
results for extremal even unimodular lattices of ranks $56$, $72$, and
$96$.  In a similar vein, Ozeki \cite{Ozeki2} showed that if $L$ is
extremal even unimodular of rank $40$, then $L$ is generated by its
vectors of norms $4$ and $6$. Here, we extend and slightly simplify Ozeki's methods, recovering Ozeki's rank-$40$ result and obtaining analogous results for extremal even unimodular lattices of ranks $80$ and $120$.

\section{Modular Forms and Theta Series}

We will use the notation $\mathcal{H} = \{z\in\mathbb{C} : \Im(z) >
0\}$ for the \textit{upper half plane} of complex numbers.
A \textit{modular form of weight $k$ for the group
  $PSL_2(\mathbb{Z})$} is a holomorphic function
$f:\mathcal{H}\to\mathbb{Z}$ which is holomorphic at $i\infty$ and
satisfies
$$
f\left(\frac{az+b}{cz+d}\right) = (cz+d)^k f(z)
$$
for all $\left(\begin{smallmatrix}a&b\\c&d\end{smallmatrix}\right)\in
PSL_2(\mathbb{Z})$. If a modular form $f$ vanishes at $z = i\infty$,
it is called a \textit{cusp form}.

Let $M_k$ and $M_k^0$ be the $\mathbb{C}$-vector spaces of modular
forms and cusp forms of weight $k$ respectively.  It is known that the
\textit{Eisenstein series}
\begin{align*}
  E_4(z) &= 1+240e^{2\pi iz} + 2160e^{4\pi iz}+6720e^{6\pi iz} + \cdots\text{ and }\\
  E_6(z) &= 1-504e^{2\pi iz} - 16632e^{4\pi iz} - 122976e^{6\pi iz} -
  \cdots,
\end{align*} which are modular forms of weights $4$ and $6$
respectively, freely generate the spaces $M_k$ in the sense that any nonzero modular form can be written uniquely as a weighted homogeneous polynomial in $E_4$ and
$E_6$. This implies that $\dim(M_k) = 0$ for $k$ odd, negative, or
$k=2$;  that $\dim(M_{2k}) = 1$ and $\dim(M_{2k}^0)=0$ for $k=0$, $2\le
k\le 5$ and $k=7$; and that multiplication by the weight-$12$
modular form $\Delta = 12^{-3}(E_4^3-E_6^2)$ defines an isomorphism
$M_{k-12}\stackrel{\sim}{\to} M_k^0$. More information on the theory
of modular forms for $PSL_2(\mathbb{Z})$ can be found in
\cite{Serre:course}.

The \textit{theta function} $\Theta_L: \mathcal{H} \to \mathbb{Z}$
associated to a lattice $L$ is defined by
$$
\Theta_L(z) = \sum_{x\in L}e^{\pi i \inprod{x}{x} z};
$$
it is a generating function encoding the norms of $L$'s vectors. For a
homogeneous \textit{harmonic polynomial} $P\in \mathbb
C[x_1,\ldots,x_n]$, i.e. a homogeneous polynomial for which
$\sum_{j=1}^n \frac{\partial^2 P}{\partial x_j^2} \equiv 0$, we define
the \textit{weighted theta series} $\Theta_{L,P}$ by
$$
\Theta_{L,P}(z) = \sum_{x\in L} P(x) e^{\pi i\inprod{x}{x} z}.
$$
As shown in \cite{Serre:course,Ebeling}, if $L$ is an even unimodular lattice of
rank $n$ then $\Theta_L$ is a modular form of weight $\frac{n}{2}$,
and if in addition $P$ is a homogeneous harmonic polynomial of degree
$d$, then $\Theta_{L,P}$ is a modular form of weight $\frac{n}{2}+d$.

\section{Main Result}\label{80!!}
We denote by $P_{d,x_0}(x)$ the ``zonal spherical harmonic
polynomial'' of degree $d$, related to the \textit{Gegenbauer
  polynomial}
by \begin{equation}\label{GEG}P_{d,x_0}(x)=G_d(\inprod{x}{x_0},
  (\inprod{x}{x}\inprod{x_0}{x_0})^{1/2}),\end{equation} 
where $G_d(\cdot, \cdot)$ is the homogeneous polynomial of degree $d$ such that
  $G_d(t,1)$ is the Gegenbauer polynomial of degree $d$ evaluated at $t$
\cite{Geg}.

We let $L$ be an extremal even unimodular lattice of rank $40r$, and
adopt the notation used by Ozeki in \cite{Ozeki2}: For an even
unimodular lattice $L$, we denote by $\Lambda_{2m}(L)$ the set of
vectors in $L$ having norm $2m$.  We denote by $\mathcal{L}_{2m}(L)$
the sublattice of $L$ generated by $\Lambda_{2m}(L)$, and similarly
denote by $\mathcal{L}_{2m_1+2m_2}(L)$ the sublattice of $L$ generated
by $\Lambda_{2m_1}(L)\cup \Lambda_{2m_2}(L)$.  

We define $a(2k,L):=|\Lambda_{2k}(L)|$.  It is clear that the theta
series $\Theta_L$ is given by $\Theta_L(z)=\sum_{k=0}^\infty
a(2k,L)e^{2k\pi i z}$.  We  note that $$4r=2\lfloor  5r/3\rfloor+2=\min\{2k>0:a(2k,L)\neq 0\}$$ is the minimal norm of vectors in $L$ and use the
notation \begin{gather*}N_j(x)=\left|\{y\in
    \Lambda_{4r}(L):\inprod{x}{y}=j\}\right|,\\M_j(x)=\left|\{y\in
    \Lambda_{4r+2}(L):\inprod{x}{y}=j\}\right|.\end{gather*}  Using the involution $y \longleftrightarrow -y$ of
$\Lambda_{m}(L)$, we see that we have $N_j(x) = N_{-j}(x)$ and $M_j(x)
= M_{-j}(x)$ for any $j\in\mathbb{R}$ and $x\in L\otimes\mathbb{R}$.

We will show the following configuration result, which directly extends Ozeki's~\cite{Ozeki2} result for extremal even
unimodular lattices of rank $40$:
\begin{thm}\label{80!}
  For $r=1,2,3$ and $L$ extremal even unimodular of rank $40r$, we
  have $L=\mathcal{L}_{4r+(4r+2)}(L).$
\end{thm}

\begin{proof}
We partition $L$ into its equivalence classes modulo
  $\mathcal{L}_{4r+(4r+2)}(L)$.  We need only show that any class $[x]\in
  L/\mathcal{L}_{4r+(4r+2)}(L)$ is represented by a vector $x_0\in [x]$
  with $\inprod{x_0}{x_0}\leq 4r+2$.

  Now, we suppose there exists some equivalence class $[x_0]\in
  L/\mathcal{L}_{4r+(4r+2)}(L)$ where $x_0\neq 0$ is a representative of
  minimal norm with $\inprod{x_0}{x_0}=2t$ for some $t\geq 2r+2$.  We
  have the inequality
  \begin{equation}\label{hop1}|\inprod{x_0}{x}|\leq 2r\text{ for all
    }x\in \Lambda_{4r}(L),\end{equation} as $x_0$ is not minimal in $L$
  whenever $\inprod{x_0}{\pm x}>2r$ since the vector $x\mp x_0$ has
  norm $$\inprod{x\mp x_0}{x\mp x_0}=\inprod{x}{x}\mp
  2\inprod{x}{x_0}+\inprod{x_0}{x_0}<\inprod{x_0}{x_0}.$$
  Similarly, we have \begin{equation}\label{hop2}|\inprod{x_0}{x}|\leq 2r+1 \text{ for all }x\in
\Lambda_{4r+2}(L).\end{equation}

  {}From (\ref{hop1}) and (\ref{hop2}), we have the
  equations \begin{eqnarray}\label{hapki1}\sum_{x\in\Lambda_{4r}(L)}\inprod{x}{x_0}^{2k}=\sum_{j=1}^{2r} 2\cdot j^{2k}\cdot N_j(x_0),\\
\label{hapki2}\sum_{x\in\Lambda_{4r+2}(L)}\inprod{x}{x_0}^{2k}=\sum_{j=1}^{2r+1}
    2\cdot j^{2k}\cdot M_j(x_0),
  \end{eqnarray}for all $k> 0$.

  We extract from the theta series $\Theta_L$ of $L$ the coefficients $a(4r,L)$ and $a(4r+2,L)$.  We observe immediately from (\ref{hapki1}) and (\ref{hapki2}) that
\begin{align}&\sum_{x\in\Lambda_{4r}(L)}\inprod{x}{x_0}^{0}=a(4r,L)\label{h-1},\\
&\sum_{x\in\Lambda_{4r+2}(L)}\inprod{x}{x_0}^{0}=a(4r+2,L)\label{h-2}.\end{align}
  Since $L$ is even unimodular of rank $40r$, we have 
  $\Theta_{L,P_{d,x_0}}\in M_{20r+d}^0$ for any $d>0$.  By comparing
  power-series coefficients, we then observe 
  \begin{align}
    \label{p1}&\Theta_{L,P_{d,x_0}}\equiv 0 \text{ for $d\in \{2,\ldots,4r-2, 4r+2\}$},\\
    \label{p2}&\Theta_{L,P_{4r,x_0}}\equiv c_1\Delta^{2r} \text{ for a constant $c_1$},\\
    \label{p3}&\Theta_{L,P_{4r+4,x_0}}\equiv c_2E_4\Delta^{2r} \text{ for a
      constant $c_2$}.
  \end{align} {}From (\ref{p1}), we obtain the equations
  \begin{align}  &\sum_{x\in\Lambda_{4r}(L)}\inprod{x}{x_0}^{2d}=a(4r,L)\frac{1\cdot 3\cdots(2d-1)}{40r\cdot
(40r+2)\cdots(40r+2d-2)}(8r)^dt^d\label{h1}\quad\text{ and}\\
&\sum_{x\in\Lambda_{4r+2}(L)}\inprod{x}{x_0}^{2d}=a(4r+2,L)\frac{1\cdot
      3\cdots(2d-1)}{40r\cdot (40r+2)\cdots(40r+2d-2)}(8r+4)^dt^d\label{h2},
  \end{align}for $d\in \{2, \ldots,4r-2, 4r+2\}$. We obtain from (\ref{p2})
   \begin{equation} \sum_{x\in\Lambda_{4r+2}(L)}P_{4r,x_0}(x) = c_{4r}
    \sum_{x\in\Lambda_{4r}(L)}P_{4r,x_0}(x),\label{h3}
  \end{equation} where $\Delta^{4r}=e^{(4r)\pi i z} + c_{4r}e^{(4r+1)\pi i  z} + O(e^{(4r+2)\pi i z})$. Similarly, (\ref{p3}) gives  \begin{equation}
    \sum_{x\in\Lambda_{4r+2}(L)}P_{4r+4,x_0}(x) = c_{4r+4} \sum_{x\in\Lambda_{4r}(L)}P_{4r+4,x_0}(x),\label{h4}
  \end{equation}where $E_4\Delta^{4r}=e^{(4r)\pi i z} + c_{4r+4}e^{(4r+1)\pi i  z} + O(e^{(4r+2)\pi i z})$.

  Combining the equations (\ref{h-1}), (\ref{h-2}), (\ref{h1}),
  (\ref{h2}), (\ref{h3}), and (\ref{h4}) with (\ref{hapki1}) and
  (\ref{hapki2}), we obtain a system of $4r+4$ homogeneous linear equations in the $4r+3$
  unknowns $$N_0(x_0),\ldots, N_{2r}(x_0), M_0(x_0),\ldots,M_{2r+1}(x_0).$$
At this stage, we diverge from our natural generalization of Ozeki's original methods and obtain the (extended) determinants of these inhomogeneous linear systems; these determinants must vanish because the system is overdetermined.

For $r=1,2,3$, these determinants are respectively{\small  \begin{gather}\label{pl1}2^{55} 3^7 5^8 7^4 11^4 13^1 19^6 23^3\cdot(t-2)\cdot t \cdot(6 t-13) \cdot\left(10 t^2-55 t+77\right),\\\label{pl2}2^{132}3^{27}5^{16}7^{10}11^{6}13^{10}23^{4}41^{8}43^{6}47^{3}\cdot (t-4)\cdot t\cdot Q_2(t),
\\2^{244} 3^{48} 5^{26} 7^{13} 11^{7} 13^{7} 17^{6} 23^{4} 31^{11} 37^{1} 59^{14} 61^{11} 67^{5} 71^{3} 73^{1}\cdot
(t-6)\cdot t\cdot Q_3(t),\label{pl3}
   \end{gather}}
where $Q_2(t)$ is the irreducible quintic \begin{gather*} 10768
  t^5-242280 t^4+2202310 t^3-10101795 t^2+23361877
  t-21771246\end{gather*} and $Q_3(t)$ is the irreducible septic \begin{gather*}19989882674056909935 t^7-892881426107875310430 t^6\\+17258039601222654151533 t^5-187053310321121904306075
   t^4\\+1227398249908229181423784 t^3-4874010945909263810320032 t^2\\+10840974078436271024624064
   t-10414527769923133690990080.\end{gather*}In each case, there are no integer solutions $t\geq 2r+2$.  However, we had assumed the existence of an equivalence class $$[x_0]\in L/\mathcal{L}_{4r+(4r+2)}(L)$$ with minimal-norm representative $x_0\neq 0$ having $\inprod{x_0}{x_0}=2t$ for integral $t\geq 2r+2$; since no such $t$ exists, all equivalence classes must be generated by vectors having norms $4r$ and $4r+2$.
\end{proof}

\section{Concluding Remarks}

A quick inspection will show that our results are the only possible immediate extensions of Ozeki's methods. In the cases $r\ge 4$, it is not possible to extract sufficiently many linear conditions by these exact techniques, as the dimensions of the relevant spaces of cusp forms grows too large.

However, using different analysis, Elkies \cite{Elky} has shown a stronger result than our Theorem \ref{80!} in the $r=3$ case: If $L$ is an extremal unimodular lattice of rank $120$ then $L = \mathcal{L}_{12}(L)$. This result for rank-$120$ lattices is analogous to Ozeki's \cite{Ozeki1, Ozeki3} results in dimensions $32$ and $48$, and to the first author's \cite{SK} results in dimensions $56$, $72$, and $96$.

\section*{Acknowledgements} The authors are especially
grateful to Professor Noam D. Elkies, under whose direction this
research was conducted.  They appreciate his instruction, support, and
commentary on drafts of this article.  The authors also thank  Denis Auroux, Katherine Paur, and Erin Schlumpf  for their translation assistance and an anonymous referee for especially helpful comments and advice.

\end{document}